\documentstyle[12pt]{article}

  \newcommand{\const}{\rm const}
  \newcommand{\Var}{\rm Var}

  \newcommand{\Dom}{\rm  Dom}

  \newcommand{\Sub}{\rm Sub}

\begin{document}

   \begin{center}

 \ {\bf  Modified Bernstein's inequality for sums of } \par

 \vspace{3mm}

 \ {\bf independent random variables.} \par

 \vspace{5mm}

 {\bf M.R.Formica, \ E.Ostrovsky, \ L.Sirota. }

\end{center}

\vspace{6mm}

 \ Universit\`{a} degli Studi di Napoli Parthenope, via Generale Parisi 13, Palazzo Pacanowsky, 80132,
Napoli, Italy. \\

e-mail: mara.formica@uniparthenope.it \\

\vspace{3mm}

 \ Israel,  Bar - Ilan University, department  of Mathematic and Statistics, 59200, \\

\vspace{3mm}

e-mails: \\
 eugostrovsky@list.ru \\
sirota3@bezeqint.net \\

\vspace{3mm}

\begin{center}

 {\bf Abstract} \\

\end{center}

\vspace{3mm}

 \ We modify the classical Bernstein's inequality for the sums of independent centered random variables (r.v.)
 in the terms of relative tails or moments.  We built also  some examples in order to show the exactness  of
offered results. \par

\vspace{4mm}

\begin{center}

\ {\it Key words and  phrases.} \par

\vspace{3mm}

\end{center}

 \hspace{3mm} Probability space, random variables (r.v.), Expectation, Variance, Banach spaces of r.v., tail characteristic,
 Banach rearrangement invariant (r.i.) functional space consisting on random variables, theory of great deviations,
 Grand Lebesgue Spaces, slowly varying function,
 Bernstein's inequality, ordinary and relative tail of distribution, independence, convexity,  Young - Fenchel
 transform, moments, norm,  norming  sum of random variables (r.v.) \par

\vspace{5mm}

 \section{Introduction. Definitions. Statement of problem.}

\vspace{3mm}

 \hspace{3mm} Let $ \ (\Omega = \{\omega\},  \ \cal{B}, \ {\bf P} ) \ $ be certain non - trivial probability space with
expectation $ {\bf E} \ $  and variance $ \ \Var. \ $ Let also $ \ \{\xi_i\}, \ i = 1,2,\ldots \ $ be a sequence of centered:
$ \ {\bf E} \xi_i =  0 \ $ independent random variables (r.v.). The classical   Bernstein's inequality, see   \cite{Bernstein},
\cite{Bernstein 2},  states that  if  this sequence satisfies the following condition

\vspace{3mm}

\begin{equation} \label{Bern condit}
\exists \kappa, \nu = \const  \in (0,\infty), \  \forall \ m = 2,3,4\ldots \hspace{3mm} \Rightarrow \  {\bf E}|\xi_i|^m \le \ \nu \ m! \ \kappa^{m - 2}/2,
\end{equation}
then

\vspace{3mm}

\begin{equation} \label{Bernst estim}
{\bf P} \left\{ \ \left| \ \sum_{i=1}^n \xi_i \ \right| > t \ \right\} \le 2 \exp \left\{ \ - \frac{t^2}{2 \nu n + 2 \kappa t}    \ \right\}, \ t \ge 0,
\end{equation}

\vspace{3mm}

 and moreover

\begin{equation} \label{Bernst maximal}
{\bf P} \left\{ \ \max_{ j = (1, 2 \ldots,n)} \  \left| \ \sum_{i=1}^j \xi_i\ \right| > t  \ \right\} \le 2 \exp \left\{ \ - \frac{t^2}{2 \nu n + 2 \kappa t}    \ \right\}, \ t \ge 0.
\end{equation}

\vspace{4mm}

 \ There are many generalizations  and applications of this result, especially in the theory of great deviations, see e.g.
\cite{Bennett 1},  \cite{Bentkus}, \cite{Kevei}, \cite{Merlev'ede}, \cite{Merlev'ede 2}, \cite{Uspensky},\cite{Utev},
\cite{Pierre Youssef} etc.\par

\vspace{4mm}

 \hspace{3mm} {\bf We offer in this short report a new approach, a new statement of problem  for uniform estimate the tail of distribution of the
 natural normalized sums of the centered independent random variables, in the terms of the so - called relative tails or moments.} \par

\vspace{4mm}

 \ We need to introduce several notations in order to formulate our statement of problem. Put $ \  \sigma_i = \sqrt{\Var ( {\xi_i} )}, \ \ $ and suppose henceforth that
 $ \ 0 < \sigma_i < \infty, \ i = 1,2,\ldots. \ $ Denote in turn

$$
S_n = \frac{\sum_{i=1}^n \xi_i}{ \ \sqrt{\sum_{i=1}^n \sigma^2_i} \ },
$$
so that $ \  \Var(S_n) = 1, \ $  a natural norming. \par

\vspace{3mm}

 \ Further, let $ \ X = (X, ||\cdot||X) \ $ be some rearrangement invariant (r.i.) Banach functional space consisting on the  numerical valued r.v. $ \ \{\eta\} \ $
having finite norms $ \ ||\eta||X < \infty,  \ $ for instance the classical Lebesgue - Riesz one $ \ L_p(\Omega, {\bf P}), \ p \ge 1, \ $ Orlicz,
Grand Lebesgue Space $  G \psi, \ $ Lorentz, Marcinkiewicz ones etc. Define for arbitrary such a space the so - called  its {\it tail characteristic} $ \ T^X(t), \ t \ge 0 \ $
 as follows:

\begin{equation} \label{Tail chact}
T^X(t) \stackrel{def}{=} \sup_{\tau \in X, \ ||\tau||X \le 1} T_{\tau}(t),
\end{equation}
where as ordinary  $ \ T_{\tau}(t) \ $  is defined as usually  the tail function for the numerical valued r.v. $ \ \tau \ $

$$
T_{\tau}(t) := {\bf P} (|\tau| \ge t), \ t \ge 0.
$$
 \vspace{3mm}

 \ As an example:  tail characteristic for the classical Lebesgue - Riesz spaces $ \  L(p) = L(p, \Omega), \ p,t \ge 1 \ $  for the sufficiently rich
 probability space. \par

 \vspace{3mm}

 \ {\bf Lemma 1.1.} \\

\begin{equation} \label{tail char for Lp}
 T^{ L(p) }(t) = t^{-p}, \ p \ge 1, t > 1.
\end{equation}

 \vspace{3mm}

  \ {\bf Proof.} The inequality  $ \  T^{ L(p) }(t) \le t^{-p}, \ p \ge 1, t > 1 \ $ follows immediately from the classical Tchebychev - Markov inequality.
  In order to ground the opposite inequality,  let us consider the following example. Let $ \ \nu \ $ be a non - negative r.v. with  the  following distribution

 $$
 {\bf P} (\nu = 0) = 1 - t^{-p}, \hspace{3mm} {\bf P} (\nu = t) = t^{-p}, \ t,p > 1.
 $$
  \ Then  $ \ {\bf E} \nu^p = 1, \ $ but  $ \ {\bf P}(\nu = t) = t^{-p}, \ $  Q.E.D. \par

\vspace{4mm}

 \hspace{3mm}  Evidently, if the non - zero r.v. $ \ \zeta \ $  belongs to the space $ \ X: \ 0 < ||\zeta||X  < \infty, \ $ then

\begin{equation} \label{General case}
T_{\zeta}(t) \le T^X(t/||\zeta||X), \ t \ge 0.
\end{equation}

 \vspace{3mm}

 \ {\bf Definition 1.1.} The r.i. space $ \ X = (X,||\cdot||X) \ $  belongs, by definition, to the class  $ \ B_2, \ $ write $ \ X \in B_2, \ $ iff
 for arbitrary independent centered r.v. $ \ \eta_1, \ \eta_2  \ $ belonging to the space $  \  X \ $ there holds

\begin{equation} \label{Pthagor inequality}
||\eta_1 + \eta_2||X \le \sqrt{(||\eta_1||X)^2 + (||\eta_2||X)^2}.
\end{equation}

 \vspace{3mm}

  \ Following, if the r.v. $ \ \{\eta_k\}, \ k = 1,2,\ldots,n \ $ are centered and commonly independent, and belonging to space $ \ X \ $ of the class $ \ B_2, \ $ then

\begin{equation} \label{key ineq}
||\Sigma_{k=1}^n \eta_k||X \le \sqrt{ \ \sum_{k=1}^n (||\eta_k||X)^2 \ }, \ n = 3,4,\ldots.
\end{equation}

\vspace{3mm}

\ For instance, the space $ \  L_2(\Omega,{\bf P}) \ $ belongs to the class $ \ B_2. \ $ \par

\vspace{3mm}

\ {\bf Definition 1.2.} The r.i. space $ \ X = (X,||\cdot||X) \ $  belongs, by definition, to the class  $ \ WB_2, \ $  (Weak $ \ B_2 \ $ space),
write $ \ X \in WB_2, \ $ iff there exists a finite constant $ \ K = K(X) \ $  such that
 for arbitrary sequence of commonly independent centered random variables $ \ \eta_1, \ \eta_2, \ldots  \ $ belonging to the space $  \  X \ $ there holds

\begin{equation} \label{key general ineq}
||\Sigma_{k=1}^n \eta_k||X \le  K(X) \cdot \sqrt{ \ \sum_{k=1}^n (||\eta_k||X)^2 \ }, \ n = 2, 3,4,\ldots.
\end{equation}

\vspace{3mm}

\ For instance, the space $ \  L_p(\Omega,{\bf P}), \ p > 2 \ $ belongs to the class $ \ WB_2. \ $ This assertion follows immediately from the
classical Rosenthal's inequality,  see \cite{Merlev'ede 3}, \ \cite{Rosenthal}. The exact value of the "constants" $ \ K(L_p), \ p > 2 \ $
is calculated in   \cite{Ibragimov Sharakhmetov}, \cite{Naimark}: \par

\begin{equation} \label{Naimark}
K(L_p) \le C[R] \cdot \frac{p}{\ln p}, \hspace{3mm} C[R] \approx 1.77638...,
\end{equation}

\vspace{3mm}

 $ \ C(R) \ $ may be named as a "Rosenthal's constant". \par

\vspace{3mm}

\ {\bf Definition 1.3.} The  {\it pair} of two r.i. spaces builded as before over source probability space,
 $ \ X = (X,||\cdot||X), \ Y = (Y, ||\cdot ||Y), \  $   is named a {\it weak pair} of the $ \ B_2 \ $ spaces,
write $ \  (X,Y) \in V\Phi_2, \ $  iff
there exists a finite constant $ \ U = U(X,Y) \ $  such that
 for arbitrary sequence of commonly independent centered r.v. $ \ \eta_1, \ \eta_2, \ldots  \ $ belonging to the space $  \  X \ $ there holds
the following estimate

\begin{equation} \label{key pair ineq}
||\Sigma_{k=1}^n \eta_k||Y \le  U(X,Y) \cdot \sqrt{ \ \sum_{k=1}^n (||\eta_k||X)^2 \ }, \ n = 2, 3,4,\ldots.
\end{equation}

 \hspace{3mm}  Evidently, if the space  $ \ X \ $ belongs to the class $ \ WB_2, \ $ then $ \ (X,X) \ \in V\Phi_2 \ $ and herewith
$ \ U(X,X) = K(X). \ $	  See also several another examples further.\\

 \vspace{4mm}

\begin{center}

 {\sc  Statement of problem.} \\

\end{center}

\vspace{4mm}

 \hspace{3mm} Suppose that the source sequence  of centered independent random variables $ \ \{\xi_i\} \ $ is such that
for some rearrangement invariant Banach functional space $ \ X = (X, ||\cdot||X) \ $ builded over $ \ (\Omega, \cal{B}, {\bf P} \ $
there holds $ \ \xi_i \in X: \ $

\begin{equation} \label{main condition}
\gamma_i := ||\xi_i/\sigma_i||X  < \infty,
\end{equation}

 \ a {\it relative  norm estimate.} \par

\vspace{4mm}

 \hspace{3mm} {\bf Our target is to obtain the uniform norm and following tail estimations for natural normed sums of these variables in the
 described above relative terms.  }\par

\vspace{4mm}

\section{Main result.}

\vspace{4mm}

 \hspace{3mm} {\bf Theorem 2.1.}  Assume that $ \ X = (X,||\cdot||X), \ Y = (Y, ||\cdot ||Y ), \  $  is the  {\it weak pair} of r.i. spaces   builded as before over
 source probability space: $ \  (X,Y) \in V\Phi_2. \ $   Let also $ \ \{\xi_i\}, \ i = 1,2,3, \ldots  \ $ be as before the sequence of independent centered non - zero
 random variables belonging to the space $ \ X.\ $ Moreover, suppose $ \ \sigma_i^2 \in (0,\infty) \ $ and that

 \begin{equation} \label{kappa restriction}
\kappa := \sup_i ||\xi_i||X/\sigma_i < \infty.
 \end{equation}

 \ Our proposition:

\vspace{3mm}

\begin{equation} \label{Prop 1}
\sup_n \ ||S_n||Y \le \kappa \cdot U(X,Y)
\end{equation}

 with correspondent uniform tail estimate  \ (\ref{General case})

\begin{equation} \label{unif tail}
\sup_n \ T_{S_n}(t) \le T^Y \ \{t/[\kappa \cdot U(X,Y)] \ \}, \ t \ge 0.
\end{equation}

\vspace{4mm}

\hspace{3mm} {\bf Proof.} Denote $ \ \xi_i = \sigma_i \ \eta_i, \ i = 1,2,\ldots;   \ $  then  $ \  ||\eta_i|| X \le \kappa, \ $

$$
S_n = \frac{\sum \xi_i}{\sqrt{\sum \sigma_i^2 }} = \frac{\sum \sigma_i \ \eta_i}{\sqrt{\sum \sigma_i^2 }};
$$

 \vspace{3mm}

$ \ \sum \stackrel{def}{=} \sum_{i=1}^n. \ $ We have

$$
||S_n||Y = ||\frac{\sum \sigma_i \ \eta_i}{\sqrt{\sum \sigma_i^2}}||Y \le
\ \frac{U(X,Y)}{\sqrt{\sum \sigma_i^2}} \cdot \sqrt{\sum \sigma^2_i \ ||\eta^2||_i X} \le
$$

$$
\frac{U(X,Y)}{\sqrt{\sum \sigma_i^2}} \cdot \sqrt{\kappa^2 \sum \sigma^2_i} = U(X,Y) \cdot \kappa,
$$
Q.E.D. \par

\vspace{4mm}

 \ {\bf Remark 2.1.} A particular case. Suppose here that the space $ \ X  = (X,||\cdot ||X ) \ $ belongs to the class $ \ WB_2. \ $
 If again $ \ \kappa \in (0,\infty), \ $ then

\vspace{3mm}

\begin{equation} \label{Prop 2}
\sup_n \ ||S_n||X \le \kappa \cdot K(X)
\end{equation}

 with correspondent uniform tail estimate

\begin{equation} \label{unif tail 2}
\sup_n \ T_{S_n}(t) \le T^X \ \{t/[\kappa \cdot K(X)] \ \}, \ t \ge 0.
\end{equation}

\vspace{4mm}

 \ {\bf Remark 2.2.} A more particular case. Suppose now that the space $ \ X  = (X,||\cdot ||X ) \ $ belongs to the class $ \ B_2. \ $
 If as before $ \ \kappa \in (0,\infty), \ $ then

\vspace{3mm}

\begin{equation} \label{Prop 3}
\sup_n \ ||S_n||X \le \kappa
\end{equation}

 with correspondent uniform tail estimate

\begin{equation} \label{unif tail 3}
\sup_n \ T_{S_n}(t) \le T^X \ \{ \ t/\kappa \}, \ t \ge 0.
\end{equation}

\vspace{4mm}

\section{Examples.   Lebesgue - Riesz and Grand Lebesgue Spaces.}

\vspace{4mm}

\begin{center}

 \hspace{3mm} {\sc A first example:  Lebesgue - Riesz spaces.} \par

\end{center}

 \vspace{3mm}

  \ Let now $ \ X = L_s(\Omega, {\bf P}), \ $  where $ \ s \in (2,\infty). \ $
  This space belongs to the class $ \ WB_2, \ $ and we deduce therefore under our notations and assumptions,
  in particular, $ \ \sigma_i \in (0,\infty), \ $

\begin{equation} \label{Leb Riesz}
\sup_n ||S_n||_s \le C_s \sup_i ||\xi_i/\sigma_i||_s, \ C_s < \infty.
\end{equation}

\vspace{3mm}

\begin{center}

\ {\sc  We bring now as an examples of these spaces for offered estimates the so - called   Grand Lebesgue Spaces (GLS). } \\

\end{center}

 \vspace{4mm}

 \ Let $\lambda_0 = \const \in (0,\infty]$ and let $ \phi = \phi(\lambda)$ be an
even strong convex function defined in the segment $  \ \lambda \in (-\lambda_0, \lambda_0) \ $ which takes only
positive values when $ \ \lambda \ne 0, \ $ twice continuously differentiable; briefly $\phi =
\phi(\lambda)$ is a Young-Orlicz function, such that

\begin{equation}\label{Young-Orlicz function}
\phi(0) = 0, \ \ \phi'(0) = 0, \ \  \phi^{''}(0) \in (0,\infty).
\end{equation}

 \ We denote the set of all these Young-Orlicz function as  $\Phi: \ \Phi = \{ \phi(\cdot)  \}. $ \\

\vspace{3mm}

{\bf Definition 3.1.}\par

\vspace{3mm}

 \hspace{3mm}  Let $\phi\in \Phi$. We say that the {\it centered} random variable $\xi$
belongs to the space $B(\phi)$  iff there exists a constant $\tau \geq 0$ such that

\begin{equation}\label{spaceB}
\forall \lambda \in (-\lambda_0, \lambda_0) \ \Rightarrow {\bf E}
\exp(\pm \lambda \ \xi) \le \exp(\phi(\lambda \ \tau)).
\end{equation}

 \ The minimal non-negative value $\tau$ satisfying ( \ref{spaceB}) for
any $\lambda \in (-\lambda_0, \ \lambda_0)$ is named $B(\phi)$-norm
of the variable $\xi$  and we write

\begin{equation}\label{Bnorm}
||\xi||_{B(\phi)} \stackrel{def}{=}\inf \{\tau \ge 0 \ : \ \forall
\lambda \in (-\lambda_0, \lambda_0) \ \Rightarrow {\bf E} \exp(\pm
\lambda \ \xi) \le \exp(\phi(\lambda \ \tau)) \} .
\end{equation}

\vspace{3mm}

 \ For instance if $\phi(\lambda)=\phi_2(\lambda) := 0.5 \ \lambda^2, \
\lambda \in \mathbf{R}$, the r.v. $\xi$ is \emph{subgaussian} and in
this case we denote the space $B(\phi_2)$ with $ \ \Sub \ $. Namely we
write $ \ \xi \in \Sub \ $ and
$$
||\xi||_{\Sub} \stackrel{def}{=} ||\xi||_{B(\phi_2)}.
$$

\vspace{4mm}

 \ It is proved in particular that $B(\phi), \ \phi  \in \Phi$, equipped with the norm
(\ref{Bnorm}) and under the ordinary algebraic operations, are
Banach rearrangement invariant  functional spaces, which are
equivalent the so-called Grand Lebesgue spaces as well as to Orlicz
exponential spaces. These spaces are very convenient for the
investigation of the r.v. having an exponential decreasing tail of
distribution; for instance, for investigation of the limit theorem,
the exponential bounds of distribution for sums of random variables,
non-asymptotical properties, problem of continuous and weak
compactness of random fields, study of Central Limit Theorem in the
Banach space, etc. \par

\vspace{3mm}

 \ Let $ \ g: R \to R \ $ be numerical valued measurable function, which can perhaps take the infinite value.
 Denote by $ \Dom[g] \ $ the domain of its finiteness:

\begin{equation}\label{Domain}
\Dom[g] := \{y, \ g(y) \in (-\infty, \ + \infty) \ \}.
\end{equation}

 \   Recall the definition $ \ g^*(u) \ $  of  the Young-Fenchel or
Legendre transform for the function $ \ g: R \to R  \ $:

\begin{equation}\label{ Definition of the Young-Fenchel transform}
g^*(u) \stackrel{def}{=} \sup_{y \in \Dom[g]} (y u - g(y)),
\end{equation}
 but we will use further the value $ \ u \ $  to be only non - negative.\par

 \ In particular, we  denote by $\nu(\cdot)$ the Young-Fenchel or
Legendre transform for the function $\phi\in \Phi$:

\begin{equation}\label{Young-Fenchel transform}
\nu(x) = \nu[\phi](x)  \stackrel{def}{=} \sup_{\lambda: |\lambda|
\le \lambda_0} (\lambda x - \phi(\lambda)) = \phi^*(x).
\end{equation}

 \ It is important to note that if the non-zero r.v. $\xi$ belongs to
the space $B(\phi)$ then

\begin{equation}\label{conditionB}
{\bf P}(\xi > x) \le \exp \left(- \nu(x/||\xi||_{B(\phi)}\right).
\end{equation}
 \ The inverse conclusion is also true up to a multiplicative constant
under  suitable conditions, see e.g. \cite{KozOsSir2017}.\par

\vspace{3mm}

\hspace{3mm} On the other words,

\begin{equation} \label{tail f Bphi}
T^{B(\phi)}(t) \le \exp \left(- \nu(t) \right), \ t \ge 0..
\end{equation}

\vspace{4mm}

 \ We recall here for reader convenience some known definitions and  facts about the so - called  Grand Lebesgue Spaces (GLS) using in this article.

   \ Let $ \ \psi = \psi(p), \ p \in [1,b) \ $    where $ \ b = \const, \ 1 <   b \le \infty \ $ be positive measurable numerical valued
    function, not necessary to be finite in every point, such that $ \ \inf_{p \in [1,b)} \psi(p) > 0. \ $   For instance

  $$
    \psi_m(p) := p^{1/m}, \ m = \const > 0, \ p \in [1,\infty)
  $$
  or

$$
   \psi^{(b; \beta)}(p) :=  (b-p)^{-\beta}, \ p \in [1,b), \ b = \const, \  1 < b < \infty; \ \beta = \const \ge 0.
$$

\vspace{4mm}

{\bf Definition 3.2.} \par

\vspace{3mm}

 \hspace{3mm}  By definition, the (Banach) Grand Lebesgue Space (GLS)    $  \ G \psi  = G\psi(b),  $
    consists on all the real (or complex) numerical valued random variable (measurable functions)
   $   \  f: \Omega \to R \ $  defined on whole our  space $ \ \Omega \ $ and having a finite norm

 \begin{equation} \label{norm psi}
    || \ f \ || = ||f||G\psi \stackrel{def}{=} \sup_{p \in [1,b)} \left[ \frac{|f|_p}{\psi(p)} \right].
 \end{equation}

 \vspace{4mm}

 \ The function $ \  \psi = \psi(p) \  $ is named as  the {\it  generating function } for this space. \par

  \ If for instance

$$
  \psi(p) = \psi^{(r)}(p) = 1, \ p = r;  \  \psi^{(r)}(p) = +\infty,   \ p \ne r,
$$
 where $ \ r = \const \in [1,\infty),  \ C/\infty := 0, \ C \in R, \ $ (an extremal case), then the correspondent
 $ \  G\psi^{(r)}(p)  \  $ space coincides  with the classical Lebesgue - Riesz space $ \ L_r = L_r(\Omega, {\bf P}). \ $ \par

\vspace{4mm}
 \ These spaces are investigated in many works, e.g. in
 \cite{Fiorenza1},   \cite{Fiorenza3}, \cite{Fiorenza4},   \cite{Iwaniec1}, \cite{Iwaniec2}, \cite{KozOs},
\cite{LiflOstSir},   \cite{Ostrovsky1}  - \cite{Ostrovsky4} etc. They are applied for example in the theory of Partial Differential Equations
\cite{Fiorenza3}, \cite{Fiorenza4}, in the theory of Probability  \cite{Ermakov},\cite{Ostrovsky3}  - \cite{Ostrovsky4}, in Statistics \cite{Ostrovsky1}, chapter 5,
theory of random fields  \cite{KozOs}, \cite{Ostrovsky4}, in the Functional Analysis \cite{Ostrovsky1}, \cite{Ostrovsky2}, \cite{Ostrovsky4} and so one. \par

 \  These spaces are rearrangement invariant (r.i.) Banach functional spaces; its fundamental function  is considered in  \cite{Ostrovsky4}. They
  not coincides  in general case with the classical  spaces: Orlicz, Lorentz, Marcinkiewicz  etc., see \cite{LiflOstSir},  \cite{Ostrovsky2}.

 \  The belonging of  some  r.v. $ \ f:  \Omega \to R \ $ to some $ \ G\psi \ $ space   is closely related with its tail behavior $ \ T_f(t) \ $
 as $ \ t\to \infty, \ $  see  \cite{KozOs}, \cite{KozOsSir2017}.    In detail, define for arbitrary generating function $ \ \psi = \psi(p) \ $

 $$
 h[\psi](p) = h(p) := p \ \ln \psi(p), \ p \in \Dom[\psi].
 $$

  \hspace{3mm}  If $ \ 0 \ne \xi \in G\psi  \ $ and if we denote $ \ ||\xi||G\psi = c \in (0,\infty), \ $ then

 \begin{equation} \label{tail Gpsi est}
 T_{\xi}(t) \le \exp ( \  - h^*( \ln (t/c)) \ ), \ t \ge ce;
 \end{equation}
and the inverse assertion holds true. Namely, it follows for arbitrary r.v. $ \ \xi \ $ from the relation  (\ref{tail Gpsi est})
that $ \ \xi \in G\psi. \ $  \par

 \vspace{4mm}

  \hspace{3mm} We conclude as before as a consequence

\begin{equation} \label{tail fun Gpsi }
 T^{G\psi} (t) \le \exp ( \  - h^*(\ln t) \ ), \ t \ge e.
 \end{equation}

\vspace{5mm}

\begin{center}

{\sc Bernstein's inequality for Grand Lebesgue Spaces. } \\

\end{center}

\vspace{5mm}

 \hspace{3mm} Let us return to the formulated above problem of modified Bernstein's inequality for the Grand Lebesgue Spaces. \par

 \vspace{3mm}

 \  It is known, see  \cite{KozOs}, \cite{Buldygin-Mushtary-Ostrovsky-Pushalsky} that if the (centered!) r.v. $\xi_i$ are
independent and subgaussian, then

\begin{equation} \label{Sums estim}
||\sum_{i=1}^n \xi_i||_{\Sub} \le \sqrt{\sum_{i=1}^n ||\xi_i||^2_{\Sub}}.
\end{equation}

\vspace{3mm}

 \ A more general statement.\par

\vspace{4mm}

\ {\bf Definition 3.3.} \ Let us impose the following important condition on the function $ \ \phi(\cdot). \ $ Indeed, we assume that
  the function  $ \ \lambda \to \phi(\sqrt{\lambda}), \ \lambda \ \in [0, \ \lambda_ 0)  \ $ is convex. Write:
  $  \ \phi(\cdot) \in \Phi(conv).  \ $ \par

\vspace{3mm}

 \ {\bf Proposition 3.1,}  see  \cite{KozOs}, \cite{KozOsSir2017}. Suppose that $ \  \phi(\cdot) \in \Phi(conv).  \ $
 Then the space $ \ B(\phi) \ $ belongs to the class $ \ B_2. \ $ \par

\vspace{4mm}

 \ Another  examples of the functions belonging to this class:

\begin{equation} \label{m L fun}
\phi_{m,L}(\lambda) \stackrel{def}{=} m^{-1} \ \lambda^m \ L(\lambda), \ \lambda \ge 1, \ m = \const \ge 2,
\end{equation}
where $ \ L = L(\lambda), \ \lambda \ge 1 \ $ is positive continuous slowly varying  at infinity function. \par

\vspace{3mm}

 \ It is known, see e.g. \ \cite{KozOsSir 3}, \ that as $ \ t \to \infty \ $

\begin{equation} \label{phi m lambda}
\phi^*_{m,L}(t) \sim g_{m,L}(t) \stackrel{def}{=} \frac{m-1}{m} \cdot t^{m/(m-1)} \cdot L^{-1/(m-1)} \left( \ t^{1/(m-1)} \ \right).
\end{equation}

\vspace{3mm}

 \hspace{3mm} The inverse proposition is also true. Indeed, suppose that some r.v. $ \ \eta \ $ is such that
$$
 T_{\eta}(t) \le c_1  \ \exp \left(- g_{m,L}(c_2 \ t) \ \right), \ t \ge 1, \
$$

then

$$
 \eta \in B_{\phi_{m,L}}: \hspace{3mm} ||\eta||B_{\phi_{m,L}} \le c_3(c_1,c_2,m,L) < \infty.
$$

\vspace{3mm}

 \ As a consequence:\\

\vspace{4mm}

 \hspace{3mm} {\bf Proposition 3.2.} It follows immediately from Remark 2.2 that if  $ \  \phi(\cdot) \in \Phi(conv)  \ $
 and

$$
\sup_i ||\xi_i/\sigma_i||B\phi = \nu \in (0,\infty),
$$
 then also

\vspace{3mm}

\begin{equation} \label{Prop 3.2}
\sup_n \ ||S_n||B\phi \le \nu
\end{equation}

 with correspondent uniform tail estimate

\begin{equation} \label{unif tail 3}
\sup_n \ T_{S_n}(t) \le T^{B\phi} \ \{ t/\nu \}, \ t \ge 0.
\end{equation}

 \vspace{4mm}

\begin{center}

 \ {\sc  The general case of Grand Lebesgue Spaces.} \\

\end{center}

\vspace{3mm}

 \hspace{3mm} Suppose now that there exists a $ \ \psi = \psi(p), \ 1 \le p < b, \ b = \const \in (1,\infty] - \ $ function
such that

$$
\sup_i ||\xi_i/\sigma_i||G\psi < \infty.
$$
 \ This generating function $ \ \psi(\cdot) \ $ may be constructively defined as follows

$$
\psi(p) \stackrel{def}{=} \sup_i ||\xi_i/\sigma_i||_p,
$$
natural way; if of course it is finite at last for certain value $ \ p > 2. \ $ \par
 \ Put also

$$
\tilde{\psi}(p) \stackrel{def}{=} C_R \ \frac{p}{\ln p} \ \psi(p), \ p \in [2,b),
$$
where $ \ b = \sup\{p, \ \psi(p) < \infty \ \}; \ $  the case when $ \ b = \infty \ $ can not be excluded.\par

\vspace{3mm}

 \ We conclude  by virtue of Rosenthal's inequality:\par

\vspace{3mm}

 \ {\bf Proposition 3.3.}

\vspace{3mm}

\begin{equation} \label{Gpsi estimat}
\sup_n ||S_n||G\tilde{\psi} \le  \sup_i ||\xi_i/\sigma_i||G\psi,
\end{equation}
with correspondent tail estimation.\par

\vspace{4mm}

\section{Non - improvability of our estimations.}

\vspace{4mm}

 \hspace{3mm} It is no hard to obtain the {\it lower estimates} for introduced here tail function $ \ T_{S_n}(t), \ t \ge 1.  \ $
Namely, note that when $ \ \sigma_i = 1, \ i = 1,2,\ldots; \ $ then the our generalized Bernstein's sum $ \ S_n \ $  coincides quite
with the classical ones

$$
S_n = n^{-1/2} \sum_{i=1}^n \xi_i,
$$
for which the lover tails estimations are well known, see
\cite{Ermakov}, \cite{Figiel}, \cite{KozOs}, \cite{Ostrovsky1}, \cite{Uspensky} etc. \par
 \ Consider for instance the following example. Let $ \ \{ \zeta_i \}, \ i = 1,2,\ldots \ $
be a sequence if independent centered identical distributed r.v. such that $ \ \sigma_i = 1 \ $
with the following tails of distributions

$$
{\bf P}(|\zeta_i| > t)  = \exp \left(- t^m \right), \ t \ge 1, \ m = \const > 1.
$$

\ Then there is a "constant" $ \  C =  C(m) \in (0,\infty) \ $ such that

\begin{equation} \label{upper m estim}
 \sup_n \ T_{S_n}(t) \le \exp \left(-  C(m) \ t^{\min(m,2)} \right), \ t \ge 1,
\end{equation}
and the power $ \ \min(m,2) \ $ is essentially non - improvable. Therefore, it is true also
for the Bernstein's statement of problem. \par

\vspace{4mm}

\section{Multivariate generalization.}

\vspace{4mm}

 \hspace{3mm} We consider in this section the case when the centered independent r.v. $ \  \{\xi_i\} \ $ taking value in
 certain separable Banach space, may be finite dimensional $ \ R^d. \ $ The Bernstein's inequality in the classical
 statement of this problem is investigated in many works: \cite{Ostrovsky5}, \cite{Prokhorov Yu}, \cite{Prokhorov A},
\cite{Yurinskii 1} - \cite{Yurinskii 4} and so one. \par

 \hspace{3mm} We introduce now the variables

 $$
 \beta_i := ||\xi_i||X, \hspace{3mm} V_n \stackrel{def}{=} \frac{\sum \xi_i}{\sqrt{\sum \beta_i^2}},
 $$
and impose as above  on the spaces $ \  X,Y \ $ in addition the following condition: for arbitrary centered independent r.v. $ \ \xi, \eta \ $
belonging to the space $ \ X \ $

\begin{equation} \label{beta cond}
||\xi + \eta||X  \le  \sqrt{||\xi||^2 + ||\eta||^2 }.
\end{equation}

\vspace{4mm}

 \ {\bf Proposition 5.1.} We find as before

\vspace{3mm}

\begin{equation} \label{mult result}
\sup_n ||V_n||Y \le K(X,Y) \ \sup_i ||\xi_i||X.
\end{equation}

\vspace{5mm}

\section{Concluding remark.}

 \hspace{3mm} It is interest in our opinion to deduce the {\it maximal} version of our estimations, for instance 
 (\ref{mult result}): to find the upper estimate  for the variable
 
 $$
 \tau(m)(t) \stackrel{def}{=} {\bf P} (\max_{n \in [1,2, \ldots, m]}   ||Y_n|| > t), \ t > 1,
 $$
in the spirit in particular of the   \cite{Szewczak}. \par

\vspace{6mm}

\vspace{0.5cm} \emph{Acknowledgement.} {\footnotesize The first
author has been partially supported by the Gruppo Nazionale per
l'Analisi Matematica, la Probabilit\`a e le loro Applicazioni
(GNAMPA) of the Istituto Nazionale di Alta Matematica (INdAM) and by
Universit\`a degli Studi di Napoli Parthenope through the project
\lq\lq sostegno alla Ricerca individuale\rq\rq (triennio 2015 - 2017)}.\par

\end{document}